\documentclass[11pt,psamsfonts]{amsart}
\usepackage[utf8]{inputenc}

\usepackage{amsmath}
\usepackage{amsthm}
\usepackage{amssymb}
\usepackage{amscd}
\usepackage{amsfonts}
\usepackage{amsbsy}
\usepackage{graphicx}
\usepackage[dvips]{psfrag}
\usepackage{array}
\usepackage{color}
\usepackage{epsfig}
\usepackage{url}

\usepackage{overpic}
\usepackage{epstopdf}

\newcommand{\R}{\ensuremath{\mathbb{R}}}

\newcommand{\x}{\mathbf{x}}

\newcommand{\bb}{\mathbf{b}}
\newcommand{\sgn}{\mathrm{sign}}
\newcommand{\de}{\delta}

\newcommand{\tr}{\mathrm{tr}}
\newcommand*{\defeq }{\mathrel{\vcenter{\baselineskip0.5ex \lineskiplimit0pt
                     \hbox{\scriptsize.}\hbox{\scriptsize.}}}%
                     =}

\newtheorem {theorem} {Theorem}

\usepackage{mathpazo}

\begin{document}
\renewcommand{\arraystretch}{1.5}

\title[A new simple proof for Lum--Chua's conjecture]
{A new simple proof for Lum--Chua's conjecture}

\author[V. Carmona, F. Fern\'{a}ndez-S\'{a}nchez, and D. D. Novaes]
{Victoriano Carmona$^1,$ Fernando Fern\'{a}ndez-S\'{a}nchez$^2,$\\ and Douglas D. Novaes$^3$}

\address{$^1$ Dpto. Matem\'{a}tica Aplicada II \& IMUS, Universidad de Sevilla, Escuela Polit\'ecnica Superior. Calle Virgen de \'Africa, 7. 41011 Sevilla, Spain.} 
\email{vcarmona@us.es}

\address{$^2$ Dpto. Matem\'{a}tica Aplicada II \& IMUS, Universidad de Sevilla, Escuela T\'{e}cnica Superior de Ingenier\'{i}a, Camino de los Descubrimientos s/n, 41092 Sevilla, Spain.} 
\email{fefesan@us.es} 

\address{$^3$Departamento de Matem\'{a}tica, Universidade
Estadual de Campinas, Rua S\'{e}rgio Buarque de Holanda, 651, Cidade Universit\'{a}ria Zeferino Vaz, 13083-859, Campinas, SP,
Brazil.} \email{ddnovaes@unicamp.br}

\subjclass[2010]{34A26,34A36,34C05,34C25}

\keywords{Piecewise planar linear systems, Limit cycles, Lum--Chua's conjecture, Poincar\'{e} half-maps}

\maketitle

\begin{abstract}
The already proved Lum--Chua's conjecture says that a continuous planar piecewise linear differential system with two zones separated by a straight line has at most one limit cycle. In this paper, we provide a new proof by using a novel characterization for Poincar\'e half-maps in planar linear systems.
This proof is very short and straightforward, because this characterization avoids the inherent flaws of the usual methods to study piecewise linear systems (the appearance of large case-by-case analysis due to the unnecessary discrimination between the different spectra of the involved matrices, to deal with transcendental equations forced by the implicit occurrence of flight time, \ldots).
In addition, the application of the characterization allow us to prove that if a limit cycle exists, then it is hyperbolic and its stability is determined by a simple relationship between  the parameters. To the best of our knowledge, the hyperbolicity of the limit cycle and this simple expression for its stability have not been pointed out before.
\end{abstract}

\section{Introduction}
\label{intro}
The study of limit cycles in planar piecewise linear differential systems dates back to Andronov et al. \cite{AndronovEtAl66} in 1937.  Since then, these systems have received a lot of attention by the scientific community mainly because of their wide range of application in applied science as idealization of nonlinear phenomenon (see, for instance, \cite{CarmonaEtAl02,FreireEtAl98,LlibreEtAl13} and the references therein).

The following continuous planar piecewise linear differential system with two zones separated by a straight line is the simplest possible configuration for a piecewise linear differential system,
\begin{equation}\label{s1}
\dot \x =
\left\{\begin{array}{l}
A_L\x+\bb, \quad\textrm{if}\quad x_1\leqslant 0,\\
A_R\x+\bb, \quad\textrm{if}\quad x_1\geqslant 0,
\end{array}\right.
\end{equation}
where $\x=(x_1,x_2)^T\in\R^2,$ $A_{L}=(a_{ij}^{L})_{2\times 2},$ $A_{R}=(a_{ij}^{R})_{2\times 2},$ with  $a_{12}^L=a_{12}^R=a_{12}$ and $a_{22}^L=a_{22}^R=a_{22},$ and $\bb=(b_1,b_2)^T\in\R^2.$

A limit cycle, in the context of continuous planar piecewise linear systems, is a non-trivial closed orbit isolated of other closed orbits (see, for instance, Definition $1.6$ in work \cite{LumChua91}). In 1991, after computer experimentations, Lum and Chua stated the following conjecture:

\medskip

\noindent{\bf Lum--Chua's Conjecture.} (\cite{LumChua91}) {\it A continuous planar piecewise differential system with two zones separated by a straight line has at most one limit cycle. The limit cycle, if it exists, is either attracting or repelling.}

\medskip

This conjecture was proven in 1998 by Freire et al. \cite{FreireEtAl98}.  Their proof is performed on a large case-by-case study, which requires individual reasonings and techniques for every possible scenario depending on
the different combinations of configurations (center, focus, saddle, node, \ldots) given by the spectra of matrices $A_L$ and $A_R$. Moreover, their approach forces the appearance of the flight-time as a new implicit variable that hinders the study.
In 2013, Llibre et al. \cite{LlibreEtAl13} made use of Massera's approach \cite{Massera54} for proving a particular case of this conjecture, namely when the determinants 
satisfy $\det(A_L)>0$ and $\det(A_R)\geqslant0$. 

Our main goal in this paper is to provide a new and simple proof for Lum--Chua's conjecture. Our proof is based on the novel integral characterization for Poincar\'e half-maps for planar linear differential systems presented in \cite{CarmonaEtAl19}, which provides a unified way to deal with the problem, avoiding the large case-by-case study of the former proof and the unnecessary appearance of the superfluous flight-time. In addition, we prove that the limit cycle, if it exists, is hyperbolic and, consequently, either attracting or repelling.  Moreover,  the stability of the limit cycle is also given in terms of the sign of a simple relationship between the parameters (see Theorem \ref{main}). As far as we know, the hyperbolicity of the limit cycle and the easy characterization of its stability have not previously been stated in the literature.

Accordingly, Lum--Chua's Conjecture follows straightforwardly from the next theorem, which is the main result of this paper.

\begin{theorem}\label{main}
The continuous planar piecewise differential system \eqref{s1} has at most one limit cycle, which is hyperbolic if it exists. Moreover, in this case, $(a_{12}b_2-a_{22}b_1) \tr(A_{L})\neq0$ and the limit cycle is attracting $($resp. repelling$)$ provided that $(a_{12}b_2-a_{22}b_1) \tr(A_{L})<0$ $($resp. $(a_{12}b_2-a_{22}b_1) \tr(A_{L})>0),$ where $\tr$ stands for the trace of the matrix.
\end{theorem}

This paper is structured as follows. First, in Section \ref{sec:NFandIIF} we introduce the preliminary results needed to prove Theorem \ref{main}. More specifically, in Section \ref{sec:NF}, we introduce the {\it Li\'{e}nard Canonical Form} for continuous piecewise linear differential systems and, in Section \ref{sec:IIF}, we  present
a result, given in \cite{CarmonaEtAl19}, that allows to manage the graphs of Poincar\'{e} half-maps of piecewise linear differential systems as orbits of a specific cubic vector field. Section \ref{sec:proof} is completely devoted to the proof of Theorem \ref{main}.
 In order to illustrate the concepts and results shown in this paper, we include and comment, in Section \ref{sec:example}, several figures that correspond to a particular system. The choice of parameters for this system is not relevant to the ideas or conclusions developed in the paper and it is obvious that any other values of the parameters would also give scenarios consistent with the proved results. A last small section is dedicated to give some conclusions.

\section{Li\'enard Canonical form and Poincar\'{e} Half-Maps}
\label{sec:NFandIIF}
This section is devoted to introduce the preliminary results needed in the proof of Theorem \ref{main}.

\subsection{Li\'{e}nard Canonical Form} \label{sec:NF}
One can readily see that the assumption $a_{12}\neq0$ is a necessary condition for the existence of periodic solutions of system \eqref{s1}. In this case, from \cite{CarmonaEtAl02}, the linear change of variables $(x,y)=(x_1,a_{22} x_1-a_{12} x_2-b_1)$ transforms system \eqref{s1} into the following Li\'{e}nard canonical form,
\begin{equation}\label{cf}
(S_L) \left\{\begin{array}{l}
\dot x= T_L x-y\\
\dot y= D_L x-a
\end{array}\right.\quad \text{for}\quad x< 0,
\quad 
(S_R) \left\{\begin{array}{l}
\dot x= T_R x-y\\
\dot y= D_R x-a
\end{array}\right.\quad \text{for}\quad x\geqslant 0,
\end{equation}
where $T_{L}=\tr(A_{L}),$ $T_{R}=\tr(A_{R}),$ $D_{L}=\det(A_{L}),$ $D_{R}=\det(A_{R}),$ and  $a=a_{12}b_2-a_{22}b_1.$ Notice that any limit cycle of system \eqref{cf} is anti-clockwise oriented and crosses the switching set 
 $\Sigma=\{(x,y)\in\mathbb{R}^2: x=0\}$ twice.

\subsection{Poincar\'{e} Half-Maps }\label{sec:IIF}

Consider the following linear differential system 
\begin{equation}\label{lienard}
\left\{\begin{array}{rl}
\dot x=&T x-y,\\
\dot y=&D x-a,
\end{array}\right.
\end{equation}
and take section 
$\Sigma=\{(x,y)\in\mathbb{R}^2: x=0\}$.

Roughly speaking, the {\it Forward Poincar\'{e} Half-Map} $y_F$ of system \eqref{lienard} associated to section 
$\Sigma$ maps a value $y_0\in[0,+\infty)$ to a value $y_1\in(-\infty,0]$, where $(0,y_1)$
is the first return to $\Sigma$ of the forward trajectory of system \eqref{lienard} starting at $(0,y_0).$ Analogously, we define the {\it Backward Poincar\'{e} Half-Map} $y_B$ for the backward trajectory. For a rigorous definition of both Poincar\'{e} half-maps, their domains and ranges, see \cite{CarmonaEtAl19}.
 
In the next section we see that $a\neq0$ is a necessary condition for the existence of limit cycles for piecewise linear system \eqref{cf}. 
Under this condition, in \cite{CarmonaEtAl19}, it is proved that the graphs of the Poincar\'{e} half-map $y_1=y_F(y_0)$ and $y_1=y_B(y_0)$  are orbits of the following cubic vector field
\begin{equation}\label{dy1}
X(y_0,y_1)=-\big(y_1 \big(Dy_0^2-aTy_0+a^2\big),y_0 \big(Dy_1^2-aTy_1+a^2\big)\big),
\end{equation}
with the same orientation as $y_0>0$ increases. 
Moreover, it is proved that the existence of the forward (resp. backward) Poincar\'{e} half-map for any value $y_0$ implies that the quadratic function $W(y)=Dy^2-aTy+a^2$ is strictly positive in the interval $[y_F(y_0),y_0]$ (resp. $[y_B(y_0),y_0]$). In particular, it is strictly positive for the domains and ranges of $y_F$ and $y_B$. 

\section{Proof of Theorem \ref{main}}\label{sec:proof}
As discussed in Section \ref{sec:NF}, system \eqref{s1} is transformed into system \eqref{cf} through a linear change of variables, provided that $a_{12}\neq0,$ which is a necessary condition for the existence of periodic solutions.

Since systems $(S_L)$ and $(S_R)$, given in  \eqref{cf}, are linear, limit cycles of system \eqref{cf} cannot be completely contained in the closed half-planes $\{(x,y)\in\mathbb{R}^2: x\leqslant 0\}$ or $\{(x,y)\in\mathbb{R}^2: x\geqslant 0\}$. Then, it is a simple consequence of the Green's Theorem that no periodic orbit exists when $T_L T_R>0$ (see \cite{FreireEtAl98}). It is also easy to see that no limit cycle can exist when either the system is homogeneous, i.e. $a=0,$ or $T_L T_R=0.$ Thus, for the sake of our interest, it is sufficient to assume that $a\neq0$ and $T_L T_R<0.$

Let $y_{L}$ (resp. $y_{R}$) be the forward (resp. backward) Poincar\'{e} half-map of planar system $(S_L)$ (resp. $(S_R)$)
associated to section $\Sigma$, 
and 
let $I_{L}\subset [0,+\infty)$ (resp. $I_{R}\subset [0,+\infty)$) be its interval of definition. 
Obviously,  no periodic solution can exist when $I_L\cap I_R=\emptyset.$ Thus, for $y_0\in I\defeq I_L\cap I_R,$ define the displacement function $\delta(y_0)= y_R(y_0)-y_L(y_0).$ Clearly, $\delta(y_0^*)=0$ if, and only if, there exists a periodic orbit passing through $(0,y_0^*)$ and $(0,y_1^*),$ $y_1^*=y_R(y_0^*)=y_L(y_0^*)$. Notice that if that periodic orbit is a limit cycle then it must be $y_0^*>0$ and $y_1^*<0$. Furthermore, it is a hyperbolic limit cycle if, and only if, $\delta'(y_0^*)\neq0.$ In this case, the limit cycle is attracting (resp. repelling) provided that $\delta'(y_0^*)<0$ (resp. $\delta'(y_0^*)>0$).

From now on, let us assume that $y_0^*\in I$, $y_0^*>0$, $\de(y_0^*)=0,$ and $y_1^*=y_R(y_0^*)=y_L(y_0^*)<0$, i.e., there exist a periodic orbit of system \eqref{cf} that crosses transversally section $\Sigma$ in two points $(0,y_0^*)$ and $(0,y_1^*)$. Next, we are going to prove that the sign of the first derivative of the displacement function $\delta$ in $y_0^*$ is given by the sign of a quadratic function at point $(y_0^*,y_1^*)$.
From \eqref{dy1}, one gets
\begin{equation}\label{diffeq}
\dfrac{d y_{L}}{d y_0}(y_0^*)=\dfrac{y_0^* W_{L}(y_1^*)}{y_1^* W_{L}(y_0^*)}
\qquad \mbox{and} \qquad
\dfrac{d y_{R}}{d y_0}(y_0^*)=\dfrac{y_0^* W_{R}(y_1^*)}{y_1^* W_{R}(y_0^*)},
\end{equation}
being
\begin{equation}\label{iif}
\begin{array}{l}
W_L(y)=D_Ly^2-aT_Ly+a^2 \qquad \mbox{and} \qquad W_R(y)=D_Ry^2-aT_Ry+a^2.
\end{array}
\end{equation}
As mentioned in Section \ref{sec:IIF}, function $W_L$ (resp. $W_R$) must be strictly positive for every interval $[y_L(y_0),y_0] $ (resp. $[y_R(y_0),y_0] $) with $y_0\in I$.

From \eqref{diffeq}, one has
$\delta'(y_0^*)=C(y_0^*,y_1^*)F(y_0^*,y_1^*),$
where functions $C$ and $F$ are real functions defined by
\[
C(y_0,y_1)=\dfrac{-y_0(y_0-y_1)}{y_1W_R(y_0)W_L(y_0)}
\]
and
\begin{equation}\label{Fe1}
F(y_0,y_1)=\dfrac{W_L(y_1)W_R(y_0)-W_L(y_0)W_R(y_1)}{y_0-y_1}.
\end{equation}
Since $W_R(y_0^*)>0$ and $W_L(y_0^*)>0$,
then the sign of $\delta'(y_0^*)$ is determined by the sign of $F.$ 
Substituting \eqref{iif} into $F(y_0,y_1),$ we get
\begin{equation*}\label{Fe2}
F(y_0,y_1)=a^3(T_L-T_R)+a(D_L T_R-D_R T_L)y_0 y_1+a^2(D_R-D_L)(y_0+y_1),
\end{equation*}
which is a quadratic function. Moreover, curve $\gamma=\{(y_0,y_1)\in\R^2:F(y_0,y_1)=0\}$ 
describes a hyperbola, possibly degenerate. Denote 
\begin{equation}\label{Q}
\begin{array}{l}
Q=\{(y_0,y_1)\in\R^2:\,
y_1\leqslant 0 \leqslant y_0\}.
\end{array}
\end{equation}
One can easily see that the curve $\gamma$ splits the set $Q\setminus\gamma$ into two disjoint connected set, $Q\setminus\gamma=R_+\cup R_-,$ where
\begin{equation}\label{Rpm}
R_{\pm}=\{(y_0,y_1)\in Q:\, \sgn(F(y_0,y_1))=\pm\sgn(aT_L)\}.
\end{equation}
Notice that $R_-$ could be the empty set, nevertheless $(0,0)\in R_+$ because $F(0,0)=a^3(T_L-T_R)$ and $T_L T_R<0$. 

Now, we are going to see that  $(y_0^*,y_1^*)\in R_+$ and, therefore, the corresponding periodic orbit is a hyperbolic limit cycle whose stability is given by $\sgn(aT_L)\neq0$. Notice that this implies the uniqueness of the limit cycle.

If $R_-=\emptyset,$ then it is trivial that $(y_0^{*},y_1^{*})\in R_+,$ so let us assume that $R_-$ is not empty.

Let us consider the graphs
$\gamma_{L}=\{(y_0,y_{L}(y_0)): \, y_0\in I_{L}\}$ and $\gamma_{R}=\{(y_0,y_{R}(y_0)): \, y_0\in I_{R}\}$, which are contained in $Q.$ Clearly, point $(y_0^*,y_1^*) \in\gamma_L\cap\gamma_R$. 
From \eqref{dy1}, $\gamma_{L}$ and $\gamma_{R}$ are orbits of the following cubic vector fields
\begin{equation}\label{X}
X_{L}(y_0,y_1)=-\big(y_1W_{L}(y_0),y_0W_{L}(y_1)\big) \ \, \mbox{and}\ \,
X_{R}(y_0,y_1)=-\big(y_1W_{R}(y_0),y_0W_{R}(y_1)\big),
\end{equation}
respectively.
 From now on, orbits $\gamma_{L}$ and $\gamma_{R}$ will be endowed with the orientation given by the corresponding vector field. These objects are illustrated in Fig.~\ref{fields} of next section for a particular example.

In the sequel, we shall study the vector fields $X_{L}$ and $X_{R}$ along the curve $\gamma$ for $(y_0,y_1)\in\textrm{int}(Q).$ From \eqref{Fe1},  if $(y_0,y_1)\in\textrm{int}(Q),$ the equation $F(y_0,y_1)=0$ is equivalent to equality $W_L(y_1)W_R(y_0)=W_L(y_0)W_R(y_1).$
Thus, substituting this last equality into $\langle\nabla F(y_0,y_1),X_{L}(y_0,y_1)\rangle$ and $\langle\nabla F(y_0,y_1),X_{R}(y_0,y_1)\rangle$ and using expression \eqref{Fe1} of $F$ we get
\begin{equation*}\label{GL}
\begin{array}{rl}
G_{L}(y_0,y_1)=&\big\langle\nabla F(y_0,y_1),X_{L}(y_0,y_1)\big\rangle\big|_{(y_0,y_1)\in \gamma
}\vspace{0.2cm}\\
=&W_{L}(y_1)a\big(T_L W_R(y_0)-T_R W_L(y_0)\big) 
\end{array}
\end{equation*}
and
\begin{equation*}\label{GR}
\begin{array}{rl}
G_{R}(y_0,y_1)=&\big\langle\nabla F(y_0,y_1),X_{R}(y_0,y_1)\big\rangle\big|_{(y_0,y_1)\in\gamma
}\vspace{0.2cm}\\
=&W_{R}(y_1)a\big(T_L W_R(y_0)-T_R W_L(y_0)\big),
\end{array}
\end{equation*}
respectively.
Since $T_L T_R<0,$ we conclude that 
\begin{equation}\label{sgnGRL}
\sgn(G_{L}(y_0,y_1))=\sgn(G_{R}(y_0,y_1))=\sgn(a T_L) \qquad \forall \ (y_0,y_1)\in\gamma.
\end{equation}
This means that each
curve $\gamma_{L}$ and $\gamma_{R}$ 
intersects $\gamma$ at most once  in $\textrm{int}(Q)$ and,  if this intersection exists, it must be crossing $\gamma$ from $R_-$ to $R_+$. This behavior is illustrated In Fig.~\ref{final} of the next section.
 
Now, since $a\neq0,$ the origin is a quadratic contact point of the continuous piecewise linear system \eqref{cf}. Thus, at least one of the Poincar\'{e} half-maps  is defined for $y_0> 0$ sufficiently small and can be continuously extended to $y_0=0$ with image $y_1=0$, see Remark 1 of  \cite{CarmonaEtAl19} (in the example considered in the next section, this is the case of the backward Poincar\'e half-map $y_R$, see Fig.~\ref{phase_planes}(b) and  Fig.~\ref{fields}(b)).
Consequently, the graph of such a Poincar\'{e} half-map contains the origin and, since $(0,0)\in R_+$, it cannot intersect the set $\gamma$ (see Fig.~\ref{final}). Hence, point $(y_0^*,y_1^*)\in \textrm{int}(Q)$  is contained in $R_+,$ that is $F(y_0^*,y_1^*)\neq0$ and $\sgn(F(y_0^*,y_1^*))=\sgn(a T_L).$ This implies that the corresponding periodic orbit is a hyperbolic limit cycle and its stability is determined by $\sgn(aT_L),$ namely it is attracting (resp. repelling) provided that $a T_L<0$ (resp. $a T_L>0$).

The proof is concluded.

\section{Illustrative example}
\label{sec:example}

Let us consider the continuous planar piecewise linear system \eqref{cf} for values
\begin{equation}\label{values}
T_L=0.4, \ T_R=-0.3, \ D_L=3, \ D_R=0.1, \ \text{and} \ a=-0.2. 
\end{equation}
Let us also consider the corresponding linear systems $(S_L)$ and $(S_R)$ involved in system \eqref{cf} for the same values given in \eqref{values}.
As it was said in last paragraph of Section \ref{intro}, the values of parameters given in \eqref{values} have been fixed just for the sake of illustrating, with nice and simple figures, some issues mentioned in the previous sections. Obviously, many other choices could also 
be taken with the same purpose.

\begin{figure}
\[
\begin{array}{cc}
\includegraphics[width=0.45\linewidth]{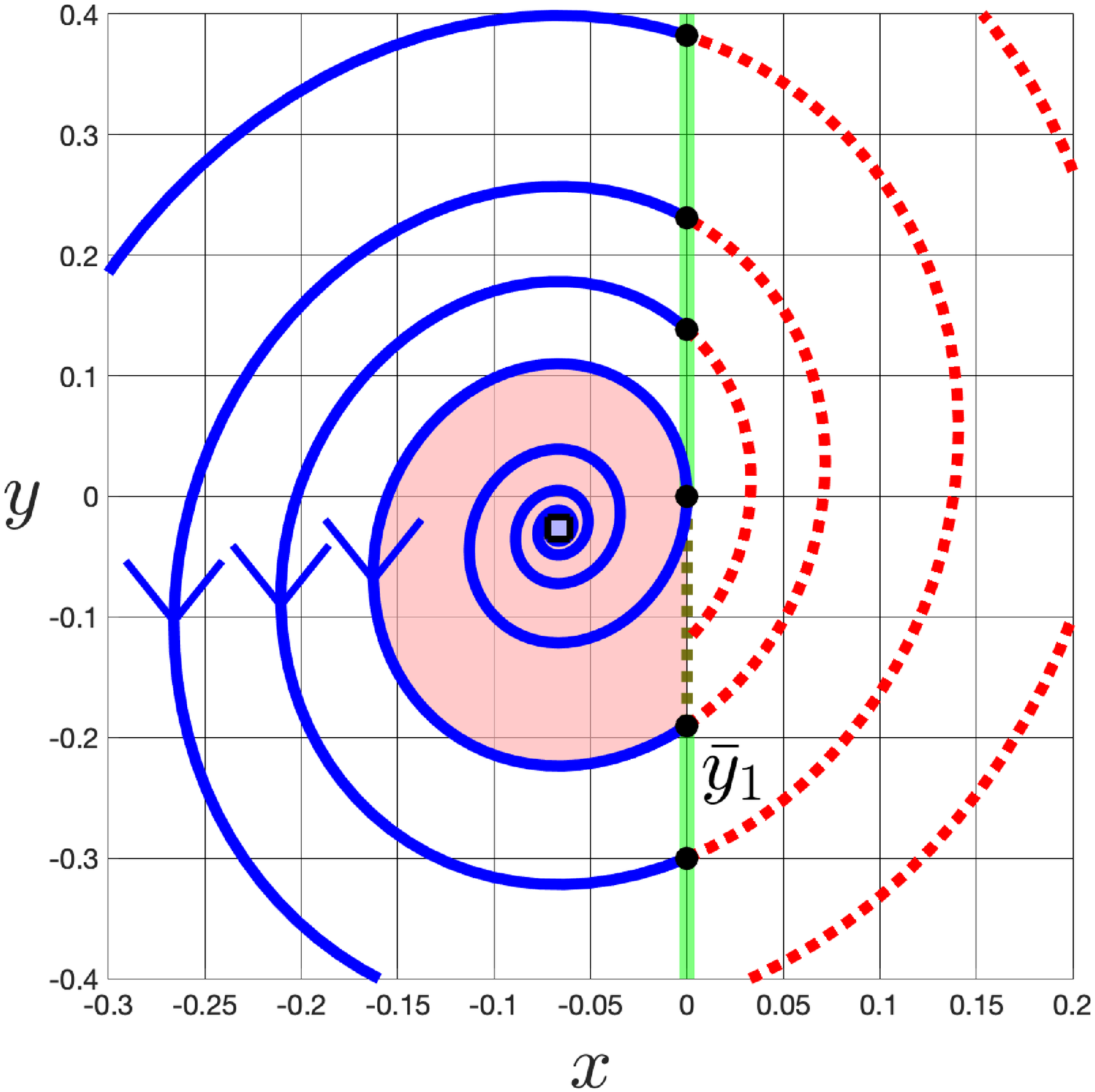}&
\includegraphics[width=0.45\linewidth]{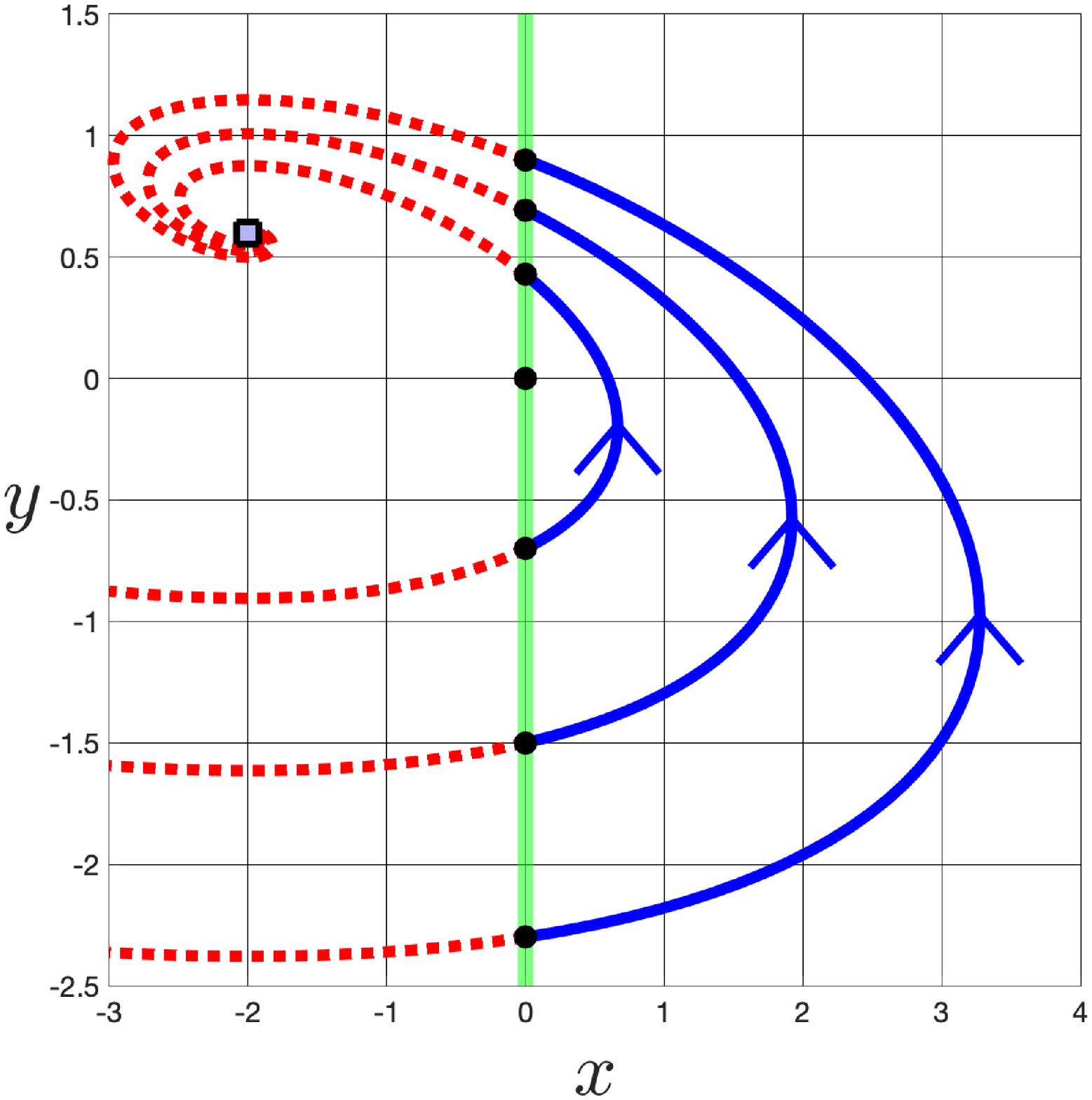}\\
\text{(a)}&\text{(b)}
\end{array}
\]
\caption{\label{phase_planes}
Phase planes of the linear systems involved in \eqref{cf} for values \eqref{values}: (a) system $(S_L)$, (b) system $(S_R)$. 
Solid curves are real pieces of orbits of piecewise system \eqref{cf}, while dotted ones stand for pieces of orbits of the linear systems $(S_L)$ and $(S_R)$ that do not correspond to orbits of system \eqref{cf}.
Note that, in (a), the interval $(\bar{y}_1,0)$, where $\bar{y}_1=y_L(0)$, is not included in the range of $y_L$. However, in (b), the backward Poincar\'e half-map $y_R$ can be continuously extended to the origin as $y_R(0)=0$, due to the quadratic contact.}
\end{figure}

The phase planes of linear systems $(S_L)$ and $(S_R)$ are shown in Fig.~\ref{phase_planes}. Solid curves are real pieces of orbits of piecewise system \eqref{cf} but dotted ones are pieces of orbits of the linear systems that do not correspond to orbits of \eqref{cf}. The orbits of piecewise system \eqref{cf} are obtained by suitable concatenations of the solid curves.

Both linear systems correspond to focus configurations. However, while the equilibrium of system $(S_L)$ is a real unstable equilibrium of system \eqref{cf}, the stable equilibrium of system $(S_R)$ is not an equilibrium of system \eqref{cf}, because it is located in the open half-plane $\{(x,y)\in\mathbb{R}: x<0\}$. Besides that,
the domain of the forward Poincar\'e half-map $y_L$ is $[0,+\infty)$
and its range is $(-\infty,\bar{y}_1]$, where $\bar{y}_1=y_L(0)$
(see Fig. \ref{phase_planes}(a)), while the domain of the
backward Poincar\'e half-map $y_R$ is $[0,+\infty)$ and its range is $(-\infty,0]$,  taking into account that it can be continuously extended to the origin as $y_R(0)=0$ (see Fig. \ref{phase_planes}(b)). Note that polynomials $W_L(y)=D_L y^2-aT_L y+a^2$ and $W_R(y)=D_R y^2-aT_R y+a^2$ are strictly positive for every $y\in\mathbb{R}$. In particular, as mentioned  in Section \ref{sec:IIF}, they are strictly positive for the domains and ranges of both Poincar\'e half-maps, as stated in \cite{CarmonaEtAl19}.

\begin{figure}
\[
\begin{array}{cc}
\includegraphics[width=0.45\linewidth]{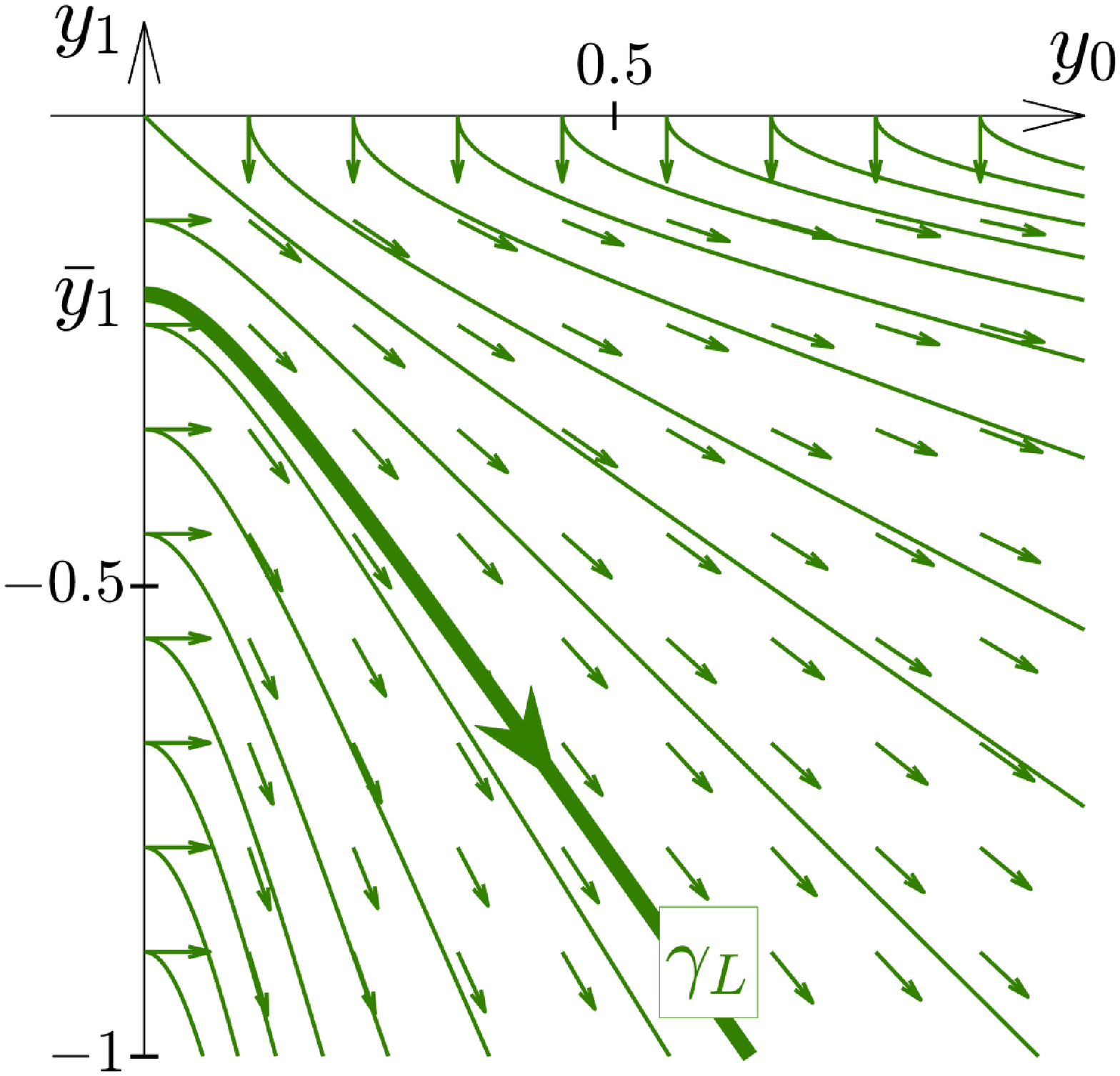}&
\includegraphics[width=0.45\linewidth]{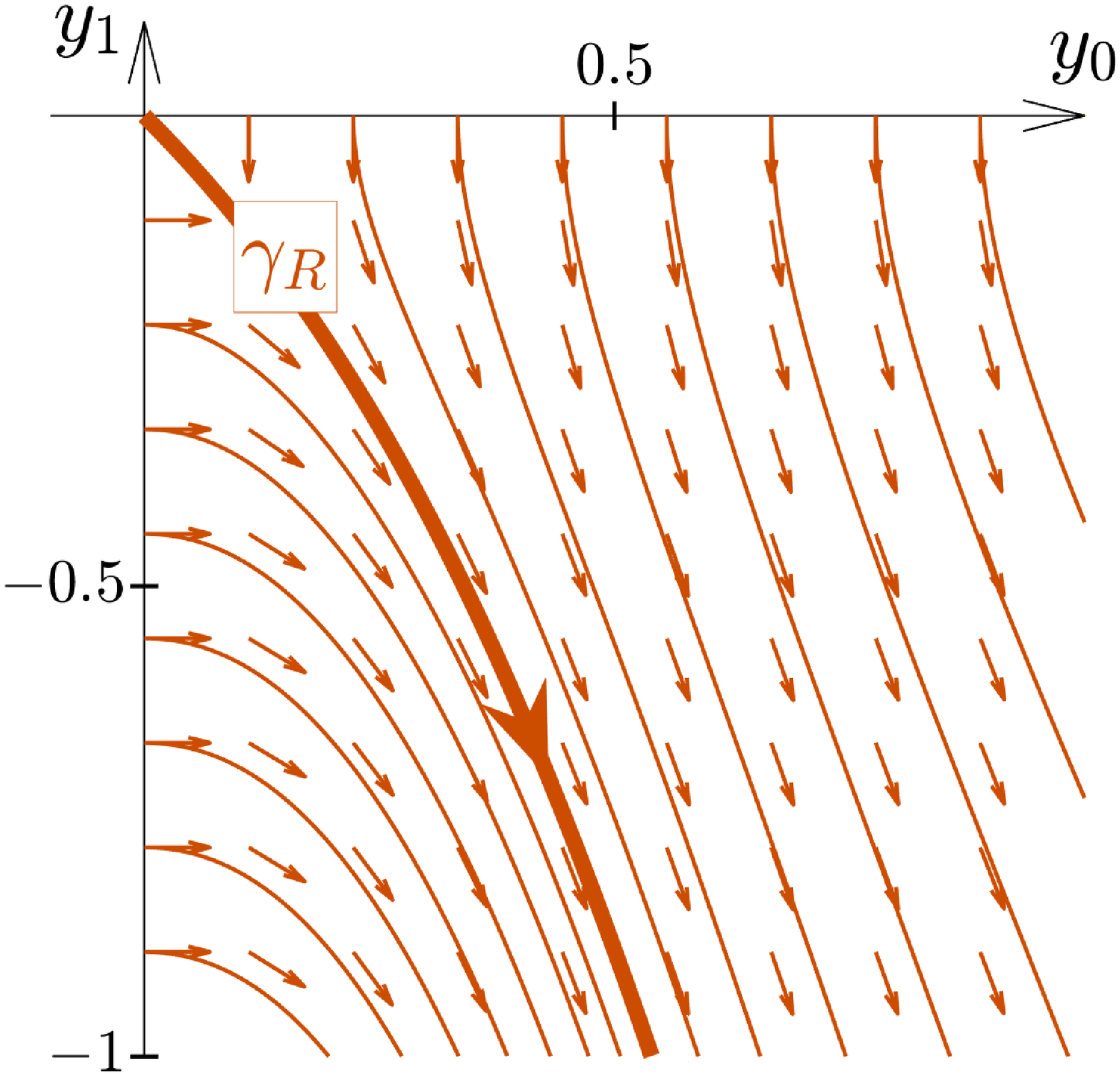}\\
\text{(a)}&\text{(b)}
\end{array}
\]
\caption{\label{fields} (a) Sketch of vector field $X_L$ given in equation \eqref{X} and several orbits of its flow (note that the vectors have been normalized for a clearer view of their directions). The thickest curve corresponds to the graph $\gamma_L$ of the forward Poincar\'e half-map $y_L$. This graph contains point $(0,\bar{y}_1)$, where $\bar{y}_1=y_L(0)$.
(b)  Idem for vector field $X_R$ given in equation \eqref{X}. Note that, due to the second order contact, the origin lies in the graph $\gamma_R$ of the backward Poincar\'e half-map $y_R$.}
\end{figure}

In Fig.~\ref{fields}, vector fields  $X_L$ and $X_R$, together with some of their orbits have been shown  in the fourth quadrant $Q$ given in \eqref{Q} (for the sake of clarity, vectors have been normalized in the figure). Among all the curves, the ones corresponding to the forward and backward Poincar\'e half-maps have been highlighted by their thickness. Taking into account the direction induced by the flow, we see that curve $\gamma_L$  begins at point $(0,\bar{y}_1)$ and $\gamma_R$  begins at the origin. This agrees with the previous descriptions of the domains and ranges of $y_L$ and $y_R$ obtained from Fig.~\ref{phase_planes}.

\begin{figure}
\includegraphics[width=0.45\linewidth]{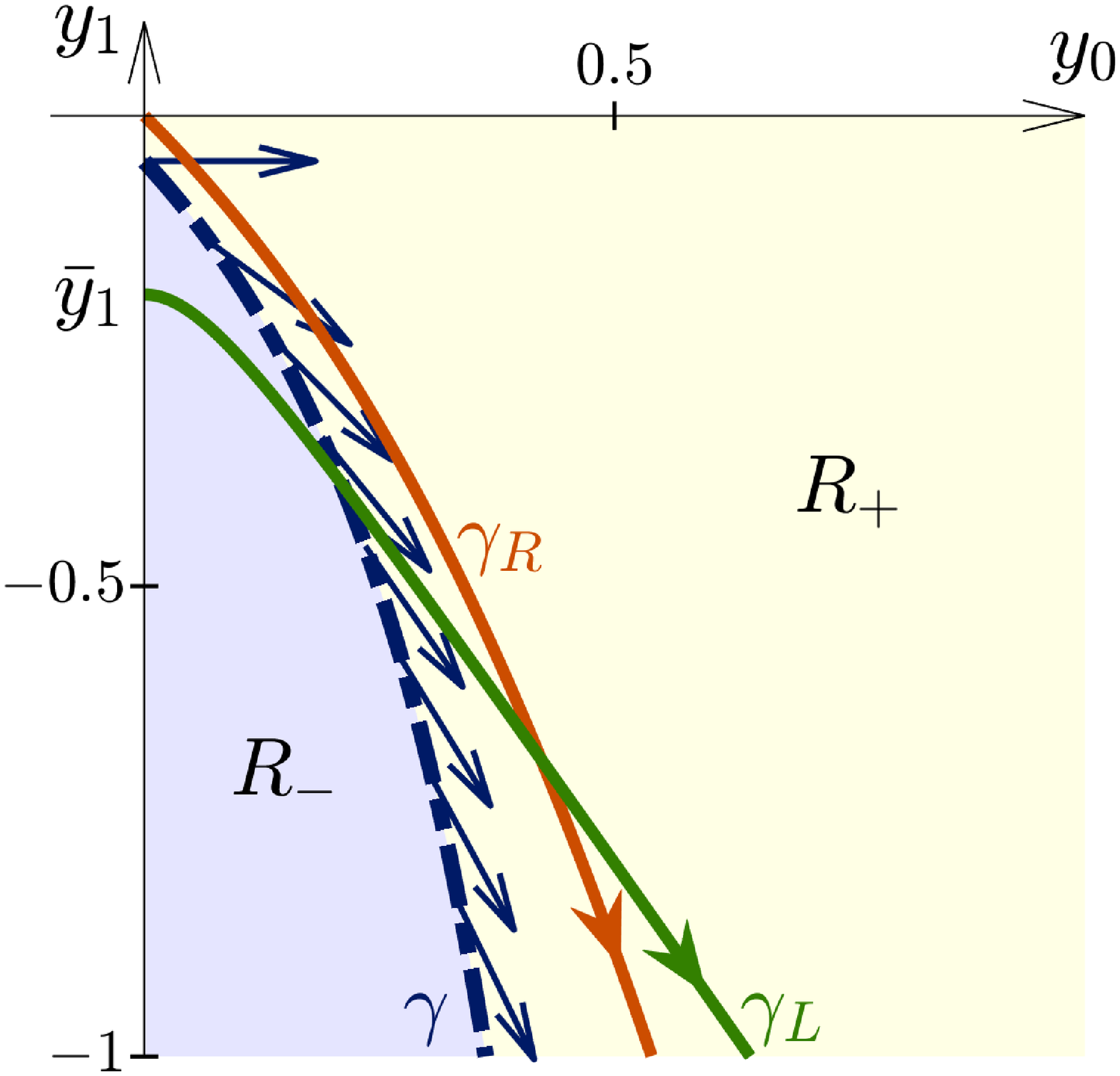}\\
\caption{\label{final} Diagram, in the fourth quadrant $Q$, of the relative position between the graphs of the Poincar\'e half-maps $y_L$ and $y_R$ (that contains the origin) and the curve $\gamma=\{(y_0,y_1)\in Q: F(y_0,y_1)=0\}$, where $F$ is defined in \eqref{Fe1}. The graph of $y_L$ contains point $(0,\bar{y}_1)$, while the origin lies in the graph of $y_R$. Vector fields $X_L$ and $X_R$, normalized like in previous figures, are sketched over the curve $\gamma$ (dashed). Actually, both normalized vector fields coincide on $\gamma$. Note that the orbits of $X_L$ and $X_R$ can cross curve $\gamma$ from $R_-$ to $R_+$. }
\end{figure}

A limit cycle for system \eqref{cf} with parameters given in \eqref{values} exists if  curves $\gamma_L$  and $\gamma_R$ have an isolated intersection point in the open set $\textrm{int}(Q)$.
Since in Fig.~\ref{final} we can observe such an intersection point, then the system has a limit cycle. 

Some other important elements for the theoretical analysis developed in the previous section have been added in Fig.~\ref{final}. In fact, the dashed curve corresponds to hyperbola $\gamma$, given by equation $F(y_0,y_1)=0$  for the values of the parameters fixed in \eqref{values}, being $F$ the function defined in \eqref{Fe1}. The vector fields $X_L$ and $X_R$ (once they have been normalized) are sketched on curve $\gamma$ in order to show that they point to $R_+$ (defined in equation \eqref{Rpm}). Notice that vector fields $X_L$ and $X_R$ have the same direction and sense over $\gamma$ (this is a direct conclusion from the equalities $F(y_0,y_1)=0$ and \eqref{sgnGRL}). Therefore, their normalizations are equal over $\gamma$.

In this example, we can check that, as was theoretically deduced in previous section, at least one of the curves of the Poincar\'e half-maps ($\gamma_R$ in this case) is always contained in region $R_+$ because it begins at the origin and cannot cross curve $\gamma$. On the contrary, curve $\gamma_L$ begins from a point of region $R_-$ and must pass through $\gamma$ in order to intersect $\gamma_R$. Since a new intersection of $\gamma_L$ and $\gamma_R$ would require that both curves pass through curve $\gamma$ once more, which is impossible, the maximum number of limit cycles must be one.

\section{Some conclusions and further directions}
\label{conclusion}
It has just been observed that the proof of Lum--Chua's conjecture can be significantly shorten and simplified, and additionally unified, by using the novel characterization of Poincar\'e half-maps given in \cite{CarmonaEtAl19}. 

This fact motivates us to analyze other problems of piecewise linear systems by using the ideas, concepts, and techniques developed in this work.

On the one hand, we would like to look for better ways to solve some previously worked problems, those involving closing equations and Poincar\'e maps, with the additional purpose of developing a common method of study. 
On the other hand, this approach allows us to consider open problems whose study with traditional techniques have turned out to be very hard to deal with.

Among the open problems, we can mention the uniqueness of limit cycles in planar piecewise linear differential systems without sliding region. This problem has been addressed by several authors, usually via large case-by-case analyses (see, for instance, \cite{FreireEtAl2013, HuanYang2013, HuanYang2014, LiLLibre2020, LlibreEtAl2008,  MedradoTorregrosa15}), some of them involving transcendental implicit equations due to the flight time. 
However, in spite of the previous efforts, as far as we know, the uniqueness of limit cycles in this case has remained an open problem. Currently, we are working in this problem and the present technique has also been proved to be very effective  \cite{CarmonaEtAl20a}.

We can also mention the problem of the existence of a finite upper bound for the maximum number of limit cycles in discontinuous planar piecewise linear differential systems. One can find some partial results in the literature for some non generic families of piecewise linear differential systems (see, for instance, \cite{FreireEtAl12,LlibreEtAl15, Novaes17}). The technique and ideas developed in the present paper allowed us to easily obtain a finite upper bound, which we are currently  improving.

\section*{Acknowledgements}

VC and FFS are partially supported by the Ministerio de Ciencia e Innovación, Plan Nacional I+D+I cofinanced with FEDER funds, in the frame of the projects MTM2014-56272-C2-1-P, MTM2015-65608-P, MTM2017-87915-C2-1-P,  PGC2018-096265-B-I00 and by the Junta de Andaluc\'{i}a grants TIC-0130, P12-FQM-1658. DDN is partially supported by Funda\c{c}\~{a}o de Amparo \`{a} Pesquisa do Estado de S\~{a}o Paulo (FAPESP) grants 2018/16430-8, 2018/ 13481-0, and 2019/10269-3, and by Conselho Nacional de Desenvolvimento Cient\'{i}fico e Tecnol\'{ o}gico CNPq grants 306649/2018-7 and 438975/ 2018-9.



\bibliography{references.bib}

\begin{thebibliography}{10}

\bibitem{AndronovEtAl66}
A.~A. Andronov, A.~A. Vitt, and S.~E. Kha\u{\i}kin.
\newblock {\em Theory of oscillators}.
\newblock Dover Publications, Inc., New York, 1987.
\newblock Translated from the Russian by F. Immirzi, Reprint of the 1966
  translation.

\bibitem{CarmonaEtAl19}
V.~Carmona and F.~Fern\'{a}ndez-S\'{a}nchez.
\newblock Integral characterization for poincar\'{e} half-maps in planar linear
  systems.
\newblock {\em arXiv:1910.13431}, 2019.

\bibitem{CarmonaEtAl20a}
V.~Carmona, F.~Fern\'{a}ndez-S\'{a}nchez, and D.~D. Novaes.
\newblock Uniqueness of limit cycles in planar piecewise linear differential
  systems without sliding region.
\newblock {\em Preprint}, 2020.

\bibitem{CarmonaEtAl02}
V.~Carmona, E.~Freire, E.~Ponce, and F.~Torres.
\newblock On simplifying and classifying piecewise-linear systems.
\newblock {\em IEEE Trans. Circuits Systems I Fund. Theory Appl.},
  49(5):609--620, 2002.

\bibitem{FreireEtAl98}
E.~Freire, E.~Ponce, F.~Rodrigo, and F.~Torres.
\newblock Bifurcation sets of continuous piecewise linear systems with two
  zones.
\newblock {\em Internat. J. Bifur. Chaos Appl. Sci. Engrg.}, 8(11):2073--2097,
  1998.

\bibitem{FreireEtAl12}
E.~Freire, E.~Ponce, and F.~Torres.
\newblock Canonical discontinuous planar piecewise linear systems.
\newblock {\em SIAM J. Appl. Dyn. Syst.}, 11(1):181--211, 2012.

\bibitem{FreireEtAl2013}
E.~Freire, E.~Ponce, and F.~Torres.
\newblock Planar filippov systems with maximal crossing set and piecewise
  linear focus dynamics.
\newblock In {\em Progress and Challenges in Dynamical Systems}, pages
  221--232. Springer Berlin Heidelberg, 2013.

\bibitem{HuanYang2013}
S.-M. Huan and X.-S. Yang.
\newblock Existence of limit cycles in general planar piecewise linear systems
  of saddle{\textendash}saddle dynamics.
\newblock {\em Nonlinear Analysis: Theory, Methods {\&} Applications},
  92:82--95, Nov. 2013.

\bibitem{HuanYang2014}
S.-M. Huan and X.-S. Yang.
\newblock On the number of limit cycles in general planar piecewise linear
  systems of node{\textendash}node types.
\newblock {\em Journal of Mathematical Analysis and Applications},
  411(1):340--353, Mar. 2014.

\bibitem{LiLLibre2020}
S.~Li and J.~Llibre.
\newblock Phase portraits of planar piecewise linear refracting systems:
  Focus-saddle case.
\newblock {\em Nonlinear Analysis: Real World Applications}, 56:103153, 2020.

\bibitem{LlibreEtAl15}
J.~Llibre, D.~D. Novaes, and M.~A. Teixeira.
\newblock Maximum number of limit cycles for certain piecewise linear dynamical
  systems.
\newblock {\em Nonlinear Dynam.}, 82(3):1159--1175, 2015.

\bibitem{LlibreEtAl13}
J.~Llibre, M.~Ord\'{o}\~{n}ez, and E.~Ponce.
\newblock On the existence and uniqueness of limit cycles in planar continuous
  piecewise linear systems without symmetry.
\newblock {\em Nonlinear Anal. Real World Appl.}, 14(5):2002--2012, 2013.

\bibitem{LlibreEtAl2008}
J.~Llibre, E.~Ponce, and F.~Torres.
\newblock On the existence and uniqueness of limit cycles in li{\'{e}}nard
  differential equations allowing discontinuities.
\newblock {\em Nonlinearity}, 21(9):2121--2142, Aug. 2008.

\bibitem{LumChua91}
R.~Lum and L.~O. Chua.
\newblock Global properties of continuous piecewise linear vector fields. part
  i: Simplest case in $\mathbb{R}^2$.
\newblock {\em International Journal of Circuit Theory and Applications},
  19(3):251--307, may 1991.

\bibitem{Massera54}
J.~L. Massera.
\newblock Sur un th\'{e}or\`eme de {G}. {S}ansone sur l'\'{e}quation di
  {L}i\'{e}nard.
\newblock {\em Boll. Un. Mat. Ital. (3)}, 9:367--369, 1954.

\bibitem{MedradoTorregrosa15}
J.~a.~C. Medrado and J.~Torregrosa.
\newblock Uniqueness of limit cycles for sewing planar piecewise linear
  systems.
\newblock {\em J. Math. Anal. Appl.}, 431(1):529--544, 2015.

\bibitem{Novaes17}
D.~D. Novaes.
\newblock Number of limit cycles for some non-generic classes of piecewise
  linear differential systems.
\newblock In {\em Extended abstracts {S}pring 2016---nonsmooth dynamics},
  volume~8 of {\em Trends Math. Res. Perspect. CRM Barc.}, pages 135--139.
  Birkh\"{a}user/Springer, Cham, 2017.

\end{thebibliography}

\bibliographystyle{abbrv}

\end{document}